\documentclass{article}
\usepackage{showlabels}
\usepackage{amsfonts}
\usepackage{amssymb}
\usepackage{amsmath}
\usepackage{graphicx}
\usepackage{amsmath,amssymb,fullpage}
\usepackage[applemac]{inputenc}
\usepackage{showlabels}

\newcommand{\dproof}{\noindent {Proof.} \quad}
\newcommand{\fproof}{\hfill $\square$ \bigskip}

\numberwithin{equation}{section}

\newenvironment{myenumerate}{%

\begin{enumerate}}{\end{enumerate}}

\def\RB{\mathbb{R}}

\def\FB{\mathbb{F}}

\def\FC{\mathcal{F}}

\def\XC{\mathcal X}

\providecommand{\U}[1]{\protect\rule{.1in}{.1in}}
\newtheorem{theorem}{Theorem}[section]

\newtheorem{proposition}[theorem]{Proposition}
\newtheorem{remark}[theorem]{Remark}

\newenvironment{proof}[1][Proof]{\noindent\textbf{#1.} }{\ \rule{0.5em}{0.5em}}

\numberwithin{equation}{section}

\def\RB{\mathbb{R}}

\def\FB{\mathbb{F}}

\def\FC{\mathcal{F}}

\def\XC{\mathcal X}

\def\1B{\text{1\!\!I}}

\def\tN{\tilde{N}}

\def\hp{\hat{p}}

\newcommand{\mr}{\mathbb{R}}
\newcommand{\cf}{{\cal F}}

\numberwithin{equation}{section}
\newtheorem{theo}{\sc Theorem}[section]
\newtheorem{ex}[theo]{\sc Example}%[section]
\newtheorem{defi}[theo]{\sc Definition}%[section]
\allowdisplaybreaks

\begin{document}

\title{Malliavin Calculus and Optimal Control of Stochastic Volterra Equations }
\author{\ Nacira Agram \thanks{Faculty of economic sciences and management, University Med
Khider, Po. Box 145, Biskra $\left(  07000\right)  $ Algeria. Email:
agramnacira@yahoo.fr} \ and Bernt \O ksendal \thanks{Corresponding author} \thanks{Dept. of Mathematics,
University of Oslo, Box 1053 Blindern, N-0316 Oslo, Norway. Email:
oksendal@math.uio.no} \thanks{Norwegian School of Economics (NHH), Helleveien 30,
N-5045 Bergen, Norway.} \ \ \ \ \ \ }

\date{\ \ 11 February 2015}

\maketitle

\begin{abstract}
Solutions of stochastic Volterra (integral) equations are not Markov processes, and therefore classical methods, like dynamic programming, cannot be used to study optimal control problems for such equations. However, we show that, by using {\em Malliavin calculus\/}, it is possible to formulate a modified functional type of {\em maximum principle\/} suitable for such systems. This principle also applies to situations where the controller has only partial information available to base her decisions upon. We present both a sufficient and a necessary maximum principle of this type, and then we use the results to study some specific examples. In particular, we solve an optimal portfolio problem in a financial market model with memory.
\end{abstract}


\textsf{MSC} \textsf{2010 Mathematics Subject Classification: 60H07, 60H20, 93E20\newline}

\section{Introduction}

Stochastic Volterra equations appear naturally in many areas of mathematics such as integral transforms, transport equations, functional differential equations and so forth, and they also appear in applications in biology, physics and finance.
For an example in economics (which also applies to population dynamics) see Example 3.4.1 in \cite{HOUZ}; and for an example stemming from Newtonian motion in a random environment, see Exercise 5.12 in \cite{O}. Stochastic Volterra equations can also be derived from stochastic \emph{delay} equations. See \cite{OZ3} and the references therein. More generally, they represent interesting models for stochastic dynamic systems with memory. For more  information on applications of Volterra integral equations, we refer to \cite{G, Belbas, OZ1} and \cite{OZ2}, the first two dealing with deterministic equations only.

In view of this, it is important to find good methods to solve optimal control problems for such equations.
In earlier papers \cite{Yong1, TQY, OZ3} and \cite{SWY}, the authors have obtained different types of maximum principles for stochastic Volterra equations.
In  \cite{Yong1} a new type of  backward stochastic Volterra integral  equations (BSVIEs), driven by Brownian motion only, is studied, and it is proved that, if a given control is optimal, then an associated BSVIE has a unique solution. (A special class of backward stochastic Volterra equations without control had earlier been studied in \cite{L}.)  An extension of the result of \cite{Yong1} to mean-field equations is obtained in \cite{SWY}.
In \cite{TQY} the same type of BSVIEs are used and  - still in the Brownian motion driven case - a necessary  maximum principle is obtained for partial information, and when the control domain is not necessarily convex.
In \cite{OZ3} a Malliavin calculus approach is used, together with a perturbation argument, to get a necessary maximum principle with partial information.

In our paper, we use Malliavin calculus to obtain both a sufficient and a necessary maximum principle for optimal control of stochastic Volterra equations with jumps and partial information. We define a Hamiltonian, which involves also the Malliavin derivatives of one of the adjoint processes. This has the advantage that the corresponding adjoint equation becomes (in some way) a standard BSDE, not a Volterra type BSVIE as in \cite{Yong1, TQY} and \cite{SWY}. On the other hand, our BSDE involves the Malliavin derivative of the adjoint process. It is interesting to note that BSDEs involving Malliavin derivatives also appear in connection with optimal control of SDEs with \emph{noisy memory}. See \cite{Dahl}.\\

Our sufficient maximum principle is new, even in the case without jumps. In our general setting also our necessary maximum principle is new. However, in the special case when the coefficients of the state equation do not depend on the state, we show that the necessary maximum principle we obtain, is equivalent to the one in \cite{OZ3}. For more general systems our maximum principle is simpler than the one in \cite{OZ3}. Moreover, there is no sufficient maximum principle in \cite{OZ3}.

In the last part of the paper, we illustrate our results by solving an optimal portfolio problem in a financial market modeled by a stochastic Volterra equation.\\

We now describe more precisely the general problem, we consider :

From now on, we let  $B(t)$ and $\tN(dt,d\zeta):= N(dt,d\zeta) -\nu(d\zeta)dt$ denote a Brownian motion and an independent  compensated Poisson random measure, respectively, on a filtered probability space
$(\Omega,   \FC, \FB:=\{\FC_t\}_{0 \leq t \leq T}, P)$ satisfying the usual conditions,  $P$ is a reference probability measure and $\nu$ is the Lvy measure of $N$. We refer to \cite{OS1} for an introduction to stochastic calculus for Lvy processes.

Let $\mathcal{A}$ be a given family of \emph{admissible} controls, required to be
$\mathcal{G}_{t}-$predictable, where $\mathbb{G=}\{\mathcal{G}_{t}\}_{t\geq0}$
is a given subfiltration of $\mathbb{F=}\{\mathcal{F}_{t}\}_{t\geq0}$, in the
sense that $\mathcal{G}_{t}\subseteq\mathcal{F}_{t}$ for all $t$.
For example, we could have%
\[
\mathcal{G}_{t}=\mathcal{F}_{(t-\delta)^{+}}\text{ \ \ \ \ (delayed
information flow).}%
\]

Suppose the state dynamics is given by a \emph{controlled stochastic Volterra
equation with jumps} of the following form:%
\begin{equation}\label{eq1.1}
\begin{array}
[c]{l}%
X(t)=X^{(u)}(t)=\xi(t)+\int_{0}^{t}b\left(  t,s,X(s),u(s)\right)  ds\\
+\int_{0}^{t}\sigma\left(  t,s,X(s),u(s)\right)  dB(s)+\int_{0}^{t}%
\int_{\mathbb{R}}\gamma\left(  t,s,X(s),u(s),\zeta\right)  \tilde{N}%
(ds,d\zeta),
\end{array}
\end{equation}
where $b(t,s,x,v)=b(t,s,x,v,\omega): [0,T]\times[0,T]\times\mathbb{R}\times \mathbb{U}\times\Omega \mapsto \mathbb{R}, \sigma(t,s,x,v)=\sigma(t,s,x,v,\omega): [0,T]\times[0,T]\times\mathbb{R}\times \mathbb{U}\times\Omega \mapsto \mathbb{R}$ and
$\gamma(t,s,x,v,\zeta)=\gamma(t,s,x,v,\zeta,\omega): [0,T]\times[0,T]\times\mathbb{R}\times \mathbb{U}\times\Omega\times\mathbb{R}_0 \mapsto \mathbb{R}$ are given functions, assumed to be $\mathbb{F}$-predictable with respect to the second variable $s$ for all $t,x,v,\zeta$ and continuously differentiable ($C^1$) with respect to the first variable $t$, with partial derivatives in $L^2([0,T]\times[0,T]\times\mathbb{R}\times\mathbb{R}\times\Omega)$ and in $L^2([0,T]\times[0,T]\times\mathbb{R}\times\mathbb{R}\times\Omega\times\nu)$, respectively. Here $\mathbb{R}_0 = \mathbb{R} - \{ 0 \} $, and $\mathbb{U}$ denotes a given open set containing all possible admissible control values $u(t,\omega)$ for
$(t,\omega) \in [0,T] \times \Omega, u \in \mathcal{A}$.

The \emph{performance functional} is given by
\begin{equation}\label{eq1.2}
J(u)=\mathbb{E}\left[
%TCIMACRO{\dint \limits_{0}^{T}}%
%BeginExpansion
{\displaystyle\int\limits_{0}^{T}}
%EndExpansion
f\left(  s,X(s),u(s)\right)  \text{ }ds+g(X(T))\right]  ; u \in \mathcal{A}.%
\end{equation}
where $f(s,x,v)=f(s,x,v,\omega): [0,T]\times\mathbb{R}\times \mathbb{U}\times\Omega \mapsto \mathbb{R}$ and $g(x)=g(x,\omega): \mathbb{R}\times\Omega \mapsto \mathbb{R}$ are given random functions, $f$ is adapted and $g$ is $\mathcal{F}_T$ - measurable and $C^1$ with respect to $x$, and we assume that $J(u)$ exists for all $u \in \mathcal{A}$.
The problem we study is the following:

\textbf{Problem:} Find $u^{\ast}\in\mathcal{A}$ such that
\begin{equation} \label{eq1.3}
J(u^{\ast})=\underset{u\in\mathcal{A}}{\sup}J(u). %
\end{equation}

Such a control $u^{\ast}$ is called an \textbf{optimal control}.

\section{A Brief Review of Malliavin Calculus for L\'{e}vy Processes}

In this section, we recall the basic definition and properties of Malliavin calculus for L\'{e}vy
processes related to this paper, for reader's convenience. A general reference for this presentation is the book \cite{DOP}. See also \cite{DMOP}, \cite{I}, \cite{N} and \cite{S}.

In view of the L\'{e}vy--It\^{o} decomposition Theorem, which states that any L\'{e}vy process
$Y(t)$ with
$$
\mathbb{E}[Y^2(t)]<\infty \quad \mbox{for all}\quad t
$$
can be written
$$
Y(t)=at+bB(t)+\int^t_0\int_{\mathbb{R}}\zeta\tilde{N}(ds,d\zeta)
$$
with constants $a$ and $b$,
we see that it suffices to deal with Malliavin calculus for $B(\cdot)$ and for
$$
\eta(\cdot):=\int_0\int_{\mathbb{R}}\zeta\tilde{N}(ds,d\zeta)
$$
separately.

\subsection{Malliavin Calculus for $B(\cdot)$}

A natural starting point is the Wiener-It\^{o} chaos expansion Theorem, which states that any
$F\in L^2({\cal F}_T,P)$ can be written
\begin{eqnarray}
F=\sum_{n=0}^{\infty}I_n(f_n)
\end{eqnarray}
for a unique sequence of symmetric deterministic functions $f_n\in L^2(\lambda^n)$,
where $\lambda$ is Lebesgue measure on $[0,T]$ and
\begin{eqnarray}
I_n(f_n)=n!\int^T_0\int^{t_n}_0\cdots\int^{t_2}_0f_n(t_1,\cdots,t_n)dB(t_1)dB(t_2)\cdots dB(t_n)
\end{eqnarray}
(the $n$-times iterated integral of $f_n$ with respect to $B(\cdot)$) for $n=1,2,\ldots$
and $I_0(f_0)=f_0$ when $f_0$ is a constant.

Moreover, we have the isometry
\begin{equation}
\mathbb{E}[F^2]=||F||^2_{L^2(p)}=\sum^\infty_{n=0}n!||f_n||^2_{L^2(\lambda^n)}.
\end{equation}

\begin{defi}[Malliavin Derivative $D_t$ with Respect to $B(\cdot)$]  %{Definition 2.1}
\hfill\break
{\rm Let
%$\mathscr{D}_{1,2}=\md^{(B)}_{1,2}$
$\mathbb{D}^{(B)}_{1,2}$ be the space of all $F\in L^2({\cf}_T,P)$
such that its chaos expansion (2.1) satisfies
\begin{eqnarray}
||F||^2_{\mathbb{D}^{(B)}_{1,2}}:=\sum^\infty_{n=1}n n!||f_n||^2_{L^2(\lambda^n)}<\infty.
\end{eqnarray}

For $F\in \mathbb{D}^{(B)}_{1,2}$ and $t\in [0,T]$, we define the {\em Malliavin derivative} (or \emph{Hida-Malliavin derivative} or \emph{the stochastic gradient}) of $F$ at $t$ (with respect to $B(\cdot)$),
$D_tF,$ by
\begin{eqnarray}\label{eq2.5a}
D_tF=\sum^\infty_{n=1}nI_{n-1}(f_n(\cdot,t)),
\end{eqnarray}
where the notation $I_{n-1}(f_n(\cdot,t))$ means that we apply the $(n-1)$-times iterated
integral to the first $n-1$ variables $t_1,\cdots, t_{n-1}$ of $f_n(t_1,t_2,\cdots,t_n)$
and keep the last variable $t_n=t$ as a parameter.}
\end{defi}

One can easily check that
\begin{eqnarray}\label{isometry}
\mathbb{E}\Big[\int^T_0(D_tF)^2dt\Big]=\sum^\infty_{n=1}n n!||f_n||^2_{L^2(\lambda^n)}=||F||^2_{\mathbb{D}^{(B)}_{1,2}},
\end{eqnarray}
so $(t,\omega)\rightarrow D_tF(\omega)$ belongs to $L^2(\lambda \times P)$.

%\underline{Example 2.2}
\begin{ex}
\rm
If $F=\int^T_0f(t)dB(t)$ with $f\in L^2(\lambda)$ deterministic, then
$$D_t F=f(t) \mbox{ for } a.a. \,t\in[0,T].$$
More generally, if $u(s)$ is Skorohod integrable, $u(s)\in \mathbb{D}_{1,2}$ for $a.a. \; s$ and $D_tu(s)$
is Skorohod integrable for $a.a. \;t$, then
\begin{equation}
D_t\Big(\int_0^Tu(s)\delta B(s)\Big)=\int_0^TD_tu(s)\delta B(s)+u(t)\;
\mbox{for a.a. $(t,\omega)$},
\end{equation}
where $\int_0^T\psi(s)\delta B(s)$ denotes the Skorohod integral of a process $\psi$
with respect to $B(\cdot)$.
%(See [N], page 35--38 for a definition of Skorohod integrals and for more details.)
\end{ex}

Some other basic properties of the Malliavin derivative $D_t$ are the following:
\begin{enumerate}
\item [(i)] {\bf Chain rule } (For a more general version see \cite{N}, page 29)\\
Suppose $F_1,\ldots,F_m\in\mathbb{D}^{(B)}_{1,2}$ and that $\psi:\mr^m\rightarrow \mr$ is $C^1$ %$e^\bot$
with bounded partial derivatives. Then, $\psi(F_1,\cdots, F_m) \in \mathbb{D}_{1,2}$ and
\begin{eqnarray}
D_t\psi(F_1,\cdots, F_m)=\sum^m_{i=1}\frac{\partial \psi}{\partial x_i}(F_1,\cdots, F_m)D_tF_i.
\end{eqnarray}

 \item [(ii)] {\bf Duality formula} \\
Suppose $u(t)$ is $\cf_t-$adapted with $\mathbb{E}[\int^T_0u^2(t)dt]<\infty$ and let $F\in \mathbb{D}^{(B)}_{1,2}$.
Then, \begin{eqnarray}\label{eq2.9a}
\mathbb{E}[F\int^T_0u(t)dB(t)]=\mathbb{E}[\int^T_0u(t)D_tFdt].
\end{eqnarray}

\item [(iii)] {\bf Malliavin derivative and adapted processes}\\
 If $\varphi $ is an $\mathbb{F}$-adapted process, then%
\[
D_{s}\varphi (t)=0\text{ for }s>t.
\]
\end{enumerate}

\begin{remark}
We put $D_{t}\varphi (t)=\underset{s\rightarrow t-}{\lim }D_{s}\varphi (t)$
(if the limit exists).
\end{remark}

\begin{remark}
It was proved in \cite{AaOPU} that one can extend the Malliavin derivative operator $D_t$ from $\mathbb{D}_{1,2}$ to all of $L^2(\mathcal{F}_T,P)$ in such a way that, also denoting the extended operator by $D_t$, for all $F \in L^2(\mathcal{F}_T,P)$ we have
\begin{equation} \label{eq2.10a}
D_tF \in (\mathcal{S})^{*} \text{ and  }
 (t,\omega) \mapsto \mathbb{E}[D_tF \mid \mathcal{F}_t] \text{  belongs to  } L^2(\lambda \times P)
\end{equation}
Here $(\mathcal{S})^{*}$ is the Hida space of stochastic distributions.\\
%and $\lambda$ denotes Lebesgue measure on $[0,T]$. \\
Moreover, the following \emph{generalized Clark-Haussmann-Ocone formula} was proved:
\begin{equation}\label{eq2.11a}
F = \mathbb{E}[F] + \int_0^T \mathbb{E}[D_tF \mid \mathcal{F}_t] dB(t)
\end{equation}
for all $F \in L^2(\mathcal{F}_T,P)$. See Theorem 3.11 in \cite{AaOPU} and also Theorem 6.35 in \cite{DOP}.\\
As also noted in \cite{Dahl} we can use this to get the following extension of the duality formula \eqref{eq2.9a}:
\end{remark}
\begin{proposition}{\bf The generalized duality formula}\\
Let $F \in L^2(\mathcal{F}_T,P)$ and let $\varphi(t,\omega) \in L^2(\lambda \times P)$ be adapted. Then
\begin{equation}
\mathbb{E}[F \int_0^T \varphi(t) dB(t)] = \mathbb{E}[ \int_0^T \mathbb{E}[D_tF \mid \mathcal{F}_t] \varphi(t) dt]
\end{equation}
\end{proposition}

\proof
By \eqref{eq2.10a} and \eqref{eq2.11a} and the It isometry we get
\begin{align}
&\mathbb{E}[F \int_0^T \varphi(t)dB(t)] = \mathbb{E}[(\mathbb{E}[F] + \int_0^T \mathbb{E}[D_tF \mid \mathcal{F}_t] dB(t))( \int_0^T \varphi(t)dB(t))]\nonumber\\
& =\mathbb{E}[ \int_0^T \mathbb{E}[D_tF \mid \mathcal{F}_t] \varphi(t) dt].
\end{align}
\fproof\\
It is this extension of the Malliavin derivative we will use from now on.
%\end{remark}

\subsection{Malliavin Calculus for $\tilde N(\cdot)$}
The construction of a stochastic derivative/Malliavin derivative in the pure jump martingale
case follows the same lines as in the Brownian motion case. In this case, the corresponding
Wiener-It\^{o} chaos expansion Theorem states that any $F\in L^2({\cf}_T,P)$ (where, in this case,
$\cf_t=\cf^{(\tilde{N})}_t$ is the $\sigma-$algebra generated by $\eta(s):=\int^s_0\int_{\mr_0}
\zeta\tilde{N}(dr,d\zeta);\; 0\leq s\leq t$) can be written as
\begin{eqnarray}\label{WIcexpan}
F=\sum^\infty_{n=0}I_n(f_n);\; f_n\in \hat{L^2}((\lambda\times \nu)^n),
\end{eqnarray}
where $\hat{L^2}((\lambda\times \nu)^n)$ is the space of functions $f_n(t_1,\zeta_1,\ldots, t_n,\zeta_n)$;
$t_i\in[0,T$], $\zeta_i\in \mathbb{R}_0$ such that 
$f_n\in L^2((\lambda\times \nu)^n)$ and $f_n$ is symmetric
with respect to the pairs of variables $(t_1,\zeta_1),\ldots,(t_n,\zeta_n).$

 It is important to note that in this case, the $n-$times iterated integral $I_n(f_n)$
 is taken with respect to $\tilde{N}(dt,d\zeta)$ and not with respect to $d\eta(t).$
Thus, we define
\begin{equation}
I_n(f_n) :=n!\int^T_0\!\!\int_{\mathbb{R}_0}\!\int^{t_n}_0\!\int_{\mathbb{R}_0}\cdots\int^{t_2}_0\!\!\int_{\mathbb{R}_0}
f_n(t_1,\zeta_1,\cdots,t_n,\zeta_n)\tilde{N}(dt_1,d\zeta_1)\cdots\tilde{N}(dt_n,d\zeta_n)
\end{equation}
for $f_n\in \hat{L^2}((\lambda\times \nu)^n).$

The It\^{o} isometry for stochastic integrals with respect to $\tilde{N}(dt,d\zeta)$
then gives the following isometry for the chaos expansion:
\begin{eqnarray}
||F||^2_{L^2(P)}=\sum^\infty_{n=0}n!||f_n||^2_{L^2((\lambda\times \nu)^n)}.
\end{eqnarray}
As in the Brownian motion case, we use the chaos expansion to define the Malliavin derivative.
Note that in this case, there are two parameters $t,\zeta,$ where $t$ represents time and $\zeta\neq 0$
represents a generic jump size.

\begin{defi}[Malliavin Derivative $D_{t,\zeta}$ with Respect to $\tilde{N}(\cdot,\cdot)$] \cite {DOP}
{\rm
Let $\mathbb{D}^{(\tilde{N})}_{1,2}$ be the space of all $F\in L^2({\cf}_T,P)$
such that its chaos expansion (\ref{WIcexpan}) satisfies
\begin{eqnarray}
||F||^2_{\mathbb{D}^{(\tilde{N})}_{1,2}}:=\sum^\infty_{n=1}n n!||f_n||^2_{L^2((\lambda\times \nu)^2)}<\infty.
\end{eqnarray}
For $F\in \mathbb{D}^{(\tilde{N})}_{1,2}$, we define the Malliavin derivative of $F$ at $(t,\zeta)$ (with respect to
$\tilde{N}(\cdot))$, $D_{t,\zeta}F, $ by
\begin{eqnarray}\label{eq2.14}
D_{t,\zeta}F:=\sum^\infty_{n=1}nI_{n-1}(f_n(\cdot,t,\zeta)),
\end{eqnarray}
where $I_{n-1}(f_n(\cdot,t,\zeta))$ means that we perform the $(n-1)-$times iterated
integral with respect to $\tilde{N}$ to the first $n-1$ variable pairs $(t_1,\zeta_1),\cdots,(t_n,\zeta_n),
$ keeping $(t_n,\zeta_n)=(t,\zeta)$ as a parameter.}
\end{defi}

In this case, we get the isometry.
\begin{eqnarray}
\mathbb{E}[\int^T_0\int_{\mathbb{R}_0}(D_{t,\zeta}F)^2 \nu (d\zeta)dt]=\sum^\infty_{n=0}n n!||f_n||^2_{L^2((\lambda\times \nu)^n)}=||F||^2_{\mathbb{D}_{1,2}^{(\tilde{N})}}.
\end{eqnarray}
(Compare with (\ref{isometry})).
%\underline{Example 2.4}

\begin{ex}
\rm
If $F=\int^T_0\int_{\mr_0}f(t,\zeta)\tilde{N}(dt,d\zeta)$
for some deterministic $f(t,\zeta)\in L^2(\lambda\times \nu)$, then
$$
D_{t,\zeta}F=f(t,\zeta) \mbox{ for } a. a. \, (t,\zeta).
$$
More  generally, if  $\psi(s,\zeta)$ is Skorohod integrable with respect to
$\tilde N(\delta s, d\zeta)$, $\psi(s,\zeta)\in \mathbb{D}_{1,2}^{(\tilde N)}$ for
$a.a.\,s,\zeta$ and $D_{t,z}\psi(s, \zeta)$ is Skorohod integrable for $a.a.\,(t,z)$, then
\begin{equation}
D_{t,z}(\int^T_0\!\int_\mr\psi(s,\zeta)\tilde{N}(\delta s,d\zeta))=\int^T_0\int_\mr D_{t,z}\psi(s,\zeta)\tilde{N}(\delta s,d\zeta)+u(t,z)\;
\mbox{ for } a.a.\, t,z,
\end{equation}
%for a.a. t,z,\\
where $\int^T_0\int_\mr\psi(s,\zeta)\tilde{N}(\delta s,d\zeta)$ denotes the
{\em Skorohod integral\/} of $\psi$ with respect to $\tilde{N}(\cdot,\cdot).$ (See \cite {DOP} for a definition of such Skorohod integrals
and for more details.)
\end{ex}

The properties of $D_{t,\zeta}$, corresponding to the properties (2.8) and (2.9) of $D_t$ are the following:
\begin{itemize}
\item [(i)] {\bf Chain rule \cite {DOP}}\\
Suppose $F_1,\cdots, F_m\in \mathbb{D}^{(\tilde{N})}_{1,2}$ and that $\phi:\mathbb{R}^m\rightarrow \mathbb{R}$
is continuous and bounded. Then, 
$\phi(F_1,\cdots,F_m)\in \mathbb{D}^{(\tilde{N})}_{1,2}$ and
\begin{equation}
D_{t,\zeta}\phi(F_1,\cdots,F_m)=\phi(F_1+D_{t,\zeta}F_1,\ldots, F_m+D_{t,\zeta}F_m)-\phi(F_1,\ldots,F_m).
\end{equation}
\item [(ii)] {\bf Duality formula \cite {DOP}}\\
%\quad
Suppose $\Psi(t,\zeta)$ is ${\cf}_t$-adapted and
$\mathbb{E}[\int^T_0\int_{\mr_0}\psi^2(t,\zeta)\nu(d\zeta)dt]<\infty$
and let $F\in \mathbb{D}_{1,2}^{(\tilde{N})}$. Then,
\begin{eqnarray}
\mathbb{E}\Big[F\int^T_0\int_{\mathbb{R}_0}\Psi(t,\zeta)\tilde{N}(dt,d\zeta)\Big]=
\mathbb{E}\Big[\int^T_0\int_{\mr_0}\Psi(t,\zeta)D_{t,\zeta}F \nu(d\zeta) dt\Big].
\end{eqnarray}
\item [(iii)]{\bf Malliavin derivative and adapted processes \cite{DOP}}\\
If $\varphi $ is an $\mathbb{F}$-adapted process, then,
\[
D_{s,\zeta }\varphi (t)=0\text{ for all }s>t.
\]
\end{itemize}
\begin{remark}
We put $D_{t,\zeta }\varphi (t)=\underset{s\rightarrow t-}{\lim }D_{s,\zeta
}\varphi (t)$ ( if the limit exists).
\end{remark}

\begin{remark}
As in Remark 2.2 we note that there is an extension of the Malliavin derivative $D_{t,\zeta}$ from  $\mathbb{D}_{1,2}^{(\tilde N)}$ to $L^2(\mathcal{F}_t \times P)$ such that the following extension of the duality theorem holds:

\begin{proposition}{\bf Generalized duality formula}\\
Suppose $\Psi(t,\zeta)$ is ${\cf}_t$-adapted and
$\mathbb{E}[\int^T_0\int_{\mr_0}\psi^2(t,\zeta)\nu(d\zeta)dt]<\infty$
and let $F\in L^2(\mathcal{F}_T \times P)$. Then,
\begin{eqnarray}
\mathbb{E}\Big[F\int^T_0\int_{\mathbb{R}_0}\Psi(t,\zeta)\tilde{N}(dt,d\zeta)\Big]=
\mathbb{E}\Big[\int^T_0\int_{\mr_0}\Psi(t,\zeta)\mathbb{E}[D_{t,\zeta}F\mid \mathcal{F}_t] \nu(d\zeta) dt\Big].
\end{eqnarray}
\end{proposition}
We refer to Theorem 13.26 in \cite{DOP}.\\
\ \\
Accordingly, note that from now on we are working with this generalized version of the Malliavin derivative.
We emphasize that this generalized Malliavin derivative $DX$ exists for all $X \in L^{2}(P)$ as an element of the Hida stochastic distribution space $(\mathcal{S})*$, and it has the property that the conditional expectation $\mathbb{E}[DX | \mathcal{F}_t]$ belongs to $L^{2}(\lambda \times P)$, where $\lambda$ is Lebesgue measure on $[0,T]$. Therefore, using this generalized Malliavin derivative, combined with conditional expectation, no assumptions on Malliavin differentiability in the classical sense are needed; we can work on the whole space of random variables in $L^{2}(P)$.
\end{remark}

\section{A Sufficient Maximum Principle}

Let $\mathcal{L}$ and $\mathcal{L}_{\zeta}$ be the set of all stochastic processes with parameter space $[0,T]$ and $[0,T]\times\mathbb{R}_0$, respectively. Define the \emph{Hamiltonian functionals}

$$H_0: [0,T] \times \mathbb{R} \times \mathbb{U} \times \mathbb{R} \times \mathbb{R} \times \mathcal{R}\mapsto  \mathbb{R}$$
and
$$H_1: [0,T] \times \mathbb{R} \times \mathbb{U} \times \mathbb{R} \times \mathcal{L} \times \mathcal{L}_{\zeta} \mapsto \mathbb{R}$$
by
\begin{equation}\label{eq2.1}
H_0(t,x,v,p,q,r):= f(t,x,v) + b(t,t,x,v)p + \sigma(t,t,x,v)q+\int_{\mathbb{R}}\gamma(t,t,x,v,\zeta)r(\zeta) \nu(d\zeta)
\end{equation}
and
\begin{align}\label{eq2.2}
H_1(t,x,v,p,D_tp(\cdot), D_{t,\zeta}p(\cdot)):= & %\xi^{\prime}(t)p(t)+
\int_{t}^{T}
\dfrac{\partial b}{\partial s}(s,t,x,v)p(s)ds+ \int_{t}^{T}
\dfrac{\partial\sigma}{\partial s}(s,t,x,v)\mathbb{E}[D_tp(s)\mid \mathcal{F}_t] ds \nonumber\\
&  + \int_{t}^{T} \int_{\mathbb{R}}
\dfrac{\partial\gamma}{\partial s}(s,t,x,v,\zeta)\mathbb{E}[D_{t,\zeta}p(s)\mid\mathcal{F}_t] \nu(d\zeta) ds.
\end{align}

Here,  $\mathcal{R}$ denotes the set of all functions
$r(\cdot):\mathbb{R}\rightarrow\mathbb{R}$
 such that the last integral above converges. We may regard $x, p, q, r=r(\cdot)$ as generic values for the processes $X(t), p(t), q(t), r(t,\cdot)$, respectively (see below).

Define
\begin{equation}\label{eq3.3}
\mathcal{H}(t,x,v,p(\cdot),q(\cdot),r(\cdot)) :=H_0(t,x,v,p,q,r) + H_1(t,x,v,p(\cdot),D_tp(\cdot),D_{t,\zeta}p(\cdot)).
\end{equation}

The BSDE for the adjoint processes $p(t),q(t),r(t,\cdot)$ is defined by%
\begin{equation}\label{eq3.4}
\begin{cases}
dp(t):=- \frac{\partial \mathcal{H}}{\partial x}(t) dt +q(t) dB(t) +\int_{\mathbb{R}} r(t,\zeta) \tilde{N}(dt,d\zeta); 0 \leq t \leq T\\
p(T) := g^{\prime}(X(T)),
\end{cases}
\end{equation}
where we have used the simplified notation
\begin{equation}\label{eq3.5}
\frac{\partial \mathcal{H}}{\partial x}(t)=\frac{\partial \mathcal{H}}{\partial x}(t,X(t),u(t),p(\cdot),q(\cdot),r(\cdot)).
\end{equation}
\\

Note that from \eqref{eq1.1} we get
\begin{equation}\label{eq3.6}
\begin{array}
[c]{l}%
dX(t)=\xi^{\prime}(t)dt+b\left(  t,t,X(t),u(t)\right)  dt+\left(
%TCIMACRO{\dint _{0}^{t}}%
%BeginExpansion
{\displaystyle\int_{0}^{t}}
%EndExpansion
\dfrac{\partial b}{\partial t}\left(  t,s,X(s),u(s)\right)  ds\right)  dt\\
+\sigma\left(  t,t,X(t),u(t)\right)  dB(t)+\left(
%TCIMACRO{\dint _{0}^{t}}%
%BeginExpansion
{\displaystyle\int_{0}^{t}}
%EndExpansion
\dfrac{\partial\sigma}{\partial t}\left(  t,s,X(s),u(s)\right)  dB(s)\right)
dt\\
+%
%TCIMACRO{\dint _{\mathbb{R}}}%
%BeginExpansion
{\displaystyle\int_{\mathbb{R}}}
%EndExpansion
\gamma\left(  t,t,X(t),u(t),\zeta\right)  \tilde{N}(dt,d\zeta)+\left(
%TCIMACRO{\dint _{0}^{t}}%
%BeginExpansion
{\displaystyle\int_{0}^{t}}
%EndExpansion%
%TCIMACRO{\dint _{\mathbb{R}}}%
%BeginExpansion
{\displaystyle\int_{\mathbb{R}}}
%EndExpansion
\dfrac{\partial\gamma}{\partial t}\left(  t,s,X(s),u(s),\zeta\right)
\tilde{N}(ds,d\zeta)\right)  dt.
\end{array}
\end{equation}

\begin{theorem}
Let $\hat{u}\in\mathcal{A},$ with corresponding solutions $\hat{X}(t),\left(
\hat{p}(t),\hat{q}(t),\hat{r}(t,\cdot)\right)  $ of \eqref{eq1.1}
 and \eqref{eq3.4}, respectively.

For fixed $t, \hat{p}(\cdot),\hat{q}(\cdot),\hat{r}(\cdot)$ define
\begin{equation}
\kappa(x,v) =  \mathcal{H}%
(t, x,v,\hat{p}(\cdot),\hat{q}(\cdot),\hat{r}(\cdot))
\end{equation}
Assume the following:
\begin{itemize}
\item
(The Arrow condition) The function
\begin{equation} \label{eq3.8A}
x \rightarrow \hat{\kappa}(x) :=\underset{v\in\mathbb{U}}{\sup}\text{ } \kappa(x,v)
\end{equation}
exists and is concave.
\item
The function%
\[
x \rightarrow g(x)  \text{ is concave.}%
\]
\item
(The maximum condition)
\begin{align}
%\label{eq2.7}
&  \underset{v\in\mathbb{U}}{\sup}\text{ }\mathbb{E}\left[  \mathcal{H}%
(t, \hat{X}(t),v,\hat{p}(\cdot),\hat{q}(\cdot),\hat{r}(\cdot))\mid
{\mathcal{G}}_{t}\right] \nonumber\\
&  =\mathbb{E}\left[  \mathcal{H}(t, \hat{X}(t),\hat{u}(t),\hat
{p}(\cdot),\hat{q}(\cdot),\hat{r}(\cdot))\mid{\mathcal{G}}_{t}\right]  \forall t   \label{eq3.9}.
\end{align}
%and that%
%\begin{equation}
%\text{(Growth conditions)} \label{1.8}%
%\end{equation}
\end{itemize}
Then, $\hat{u}$ is an optimal control.
\end{theorem}

\dproof
By considering a suitable increasing family of stopping times converging to $T$, we may assume that all the local martingales appearing in the proof below are martingales. See the proof of Theorem 2.1 in \cite{OS2} for details. \\
Choose an arbitrary $u\in\mathcal{A}$ with corresponding $X(t)$ and consider
\begin{equation}
J(u)-J(\hat{u})=I_{1}+I_{2},\nonumber
\end{equation}
where
\begin{equation}\label{eq2.8}
I_{1}=\mathbb{E}\left[  \int\limits_{0}^{T}\left\{  f\left(  t\right)
-\hat{f}\left(  t\right)  \right\}  dt\text{ }\right], \quad I_{2}=\mathbb{E}%
\left[  g\left(  X(T\right)  )-g(\hat{X}\left(  T\right)  )\right]  ,
\end{equation}

where $f(t)=$ $f\left(  t,X(t),u(t)\right)  ,$ $\hat{f}\left(  t\right)
=f(t,\hat{X}(t),\hat{u}(t)).$

Using a similar notation for
\[
b(t)=b\left(  t,t,X(t),u(t)\right)  ,\hat{b}\left(  t\right)  =b(t,t,\hat
{X}(t),\hat{u}(t))\text{ etc.,}%
\]

we get
\begin{align}\label{eq2.9}
I_{1}  &  =\mathbb{E} [  \int_{0}^{T} \{ H_0(t,X(t),u(t),\hat{p}(t),\hat{q}(t),\hat{r}(t,\cdot)) \nonumber\\
&  -H_0(t,\hat{X}(t),\hat{u}(t),\hat{p}(t),\hat{q}(t),\hat{r}(t,\cdot)) \nonumber\\
&  -[ b(t)-\hat{b}(t)]\hat{p}(t) - [ \sigma(t)  -\hat{\sigma}(t)]  \hat{q}(t)-\int_{\mathbb{R}}[  \gamma (t,\zeta)
 -\hat{\gamma}(t,\zeta)]  \hat{r}
(t,\zeta)\nu(d\zeta)\} dt ].
\end{align}

Using concavity and the It\^ o Formula, we obtain%
\begin{align}\label{eq2.10}
I_{2}  &  \leq\mathbb{E}\left[  g^{\prime}(\hat{X}(T))\left(  X(T)-\hat
{X}(T)\right)  \right] \nonumber\\
&  =\mathbb{E}\left[  \hat{p}(T)\left(  X(T)-\hat{X}(T)\right)  \right]
\nonumber\\
&  =\mathbb{E}\left[  \int_{0}^{T}\left\{  \hat{p}(t)\left(  b(t)-\hat
{b}\left(  t\right)  +\int_{0}^{t}\left(  \frac{\partial b}{\partial t}\left(
t,s\right)  -\frac{\partial\hat{b}}{\partial t}\left(  t,s\right)  \right)
ds\right.  \right.  \right. \nonumber\\
&  +\int_{0}^{t}\left(  \frac{\partial\sigma}{\partial t}\left(  t,s\right)
-\frac{\partial\hat{\sigma}}{\partial t}\left(  t,s\right)  \right)
dB(s)\nonumber\\
&  \left.  +\int_{0}^{t}\int_{\mathbb{R}}\left(  \frac{\partial\gamma}{\partial t}\left(
t,s,\zeta\right)  -\frac{\partial\hat{\gamma}}{\partial t}\left(
t,s,\zeta\right)  \right)  \tilde{N}(ds,d\zeta) \right) \nonumber\\
&   - \frac{\partial \hat{\mathcal{H}}}{\partial x}(t) \left( X(t)-\hat{X}(t)\right)  \nonumber\\
&  \left.  \left.  +\hat{q}(t)\left[  \sigma\left(  t\right)  -\hat{\sigma
}\left(  t\right)  \right]  +\int_{\mathbb{R}}\hat{r}(t,\zeta)\left[
\gamma\left(  t,\zeta\right)  -\hat{\gamma}\left(  t,\zeta\right)  \right]
\nu(d\zeta)\right\}  dt\right]  .
\end{align}

By the Fubini Theorem, we get
\begin{equation}\label{eq2.11}
\int_0^T(\int_0^t \frac{\partial b}{\partial t}(t,s) ds)\hp(t)dt = \int_0^T(\int_s^T \frac{\partial b}{\partial t}(t,s) \hp(t) dt)ds
= \int_0^T(\int_t^T \frac{\partial b}{\partial s}(s,t) \hp(s) ds)dt,
\end{equation}

and similarly, by the duality theorems,
\begin{align}\label{eq2.12}
\mathbb{E}[\int_0^T(\int_0^t \frac{\partial \sigma}{\partial t}(t,s) dB(s))\hp(t)dt] = \int_0^T\mathbb{E}[\int_0^t \frac{\partial \sigma}{\partial t}(t,s) dB(s)\hp(t)]dt
=\int_0^T\mathbb{E}[\int_0^t \frac{\partial \sigma}{\partial t}(t,s) \mathbb{E}[D_s\hp(t)\mid\mathcal{F}_s] ds]dt \nonumber\\
= \int_0^T\mathbb{E}[\int_s^T \frac{\partial \sigma}{\partial t}(t,s) \mathbb{E}[D_s\hp(t)\mid\mathcal{F}_s] dt] ds
= \mathbb{E}[\int_0^T\int_t^T \frac{\partial \sigma}{\partial s}(s,t) \mathbb{E}[D_t\hp(s)\mid\mathcal{F}_t] ds dt]
\end{align}
and

\begin{align}\label{eq2.13}
\mathbb{E}[\int_0^T(\int_0^t \int_{\mathbb{R}} \frac{\partial \gamma}{\partial t}(t,s,\zeta)  \tilde{N}(ds,d\zeta))\hp(t)dt]
=\int_0^T\mathbb{E}[\int_0^t \int_{\mathbb{R}} \frac{\partial \gamma}{\partial t}(t,s,\zeta)  \tilde{N}(ds,d\zeta))\hp(t)]dt\nonumber\\
= \int_0^T\mathbb{E}[\int_0^t \int_{\mathbb{R}}\frac{\partial \gamma}{\partial t}(t,s,\zeta) \mathbb{E}[D_{s,\zeta}\hp(t)\mid\mathcal{F}_s] \nu(d\zeta)ds]dt
= \int_0^T\mathbb{E}[\int_s^T \int_{\mathbb{R}}\frac{\partial \gamma}{\partial t}(t,s,\zeta) \mathbb{E}[D_{s,\zeta}\hp(t)\mid\mathcal{F}_s] \nu(d\zeta)dt]ds\nonumber\\
= \mathbb{E}[\int_0^T\int_t^T \int_{\mathbb{R}}\frac{\partial \gamma}{\partial s}(s,t,\zeta) \mathbb{E}[D_{t,\zeta}\hp(s)\mid\mathcal{F}_t]\nu(d\zeta)dsdt]. \nonumber\\
\end{align}

Substituting \eqref{eq2.12} and \eqref{eq2.13} into \eqref{eq2.10}, we get
\begin{align}\label{eq2.16}
I_2 & \left. \leq \mathbb{E}\left[  \int_{0}^{T}\left\{
\hat{p}(t)
[b(t)-\hat
{b}\left(  t\right)]  + \hp(t)\int_{0}^{t}
(\frac{\partial b}{\partial t}\left(
t,s\right)  -\frac{\partial\hat{b}}{\partial t}\left(  t,s\right))
ds \right.  \right.  \right. \nonumber\\
&  + \int_t^T (\frac{\partial \sigma}{\partial s}(s,t)-\frac{\partial \hat{\sigma}}{\partial s}(s,t)) \mathbb{E}[D_t \hp(s)\mid\mathcal{F}_t] ds\nonumber\\
&  \left.  +\int_t^T \int_{\mathbb{R}}(\frac{\partial \gamma}{\partial s}(s,t,\zeta)-\frac{\partial \hat{\gamma}}{\partial s}(s,t,\zeta) ) \mathbb{E}[D_{t,\zeta} \hp(s)\mid\mathcal{F}_t] ds \nu(d\zeta) \right.
\nonumber\\
&   -  \frac{\partial \hat{\mathcal{H}}}{\partial x}(t) \left( X(t)-\hat{X}(t)\right)  \nonumber\\
&  \left.  \left.  +\hat{q}(t)\left[  \sigma\left(  t\right)  -\hat{\sigma
}\left(  t\right)  \right]  +\int_{\mathbb{R}}\hat{r}(t,\zeta)\left[
\gamma\left(  t,\zeta\right)  -\hat{\gamma}\left(  t,\zeta\right)  \right]
\nu(d\zeta)\right\}
dt\right]  .
\end{align}

Adding \eqref{eq2.9} and \eqref{eq2.16}, we get
\begin{align}
J(u)-J(\hat{u})  &  =I_{1}+I_{2}\\
&  \leq\mathbb{E}\left[  \int_{0}^{T}\left\{  \mathcal{H}(t)-\widehat
{\mathcal{H}}(t)- \frac{\partial \hat{\mathcal{H}}}{\partial x}(t) \left(
X(t)-\hat{X}(t)\right)  \right\}  dt\right]. \nonumber
\end{align}

Finally, to prove that
\begin{equation}
 \mathcal{H}(t)-\widehat
{\mathcal{H}}(t)- \frac{\partial \hat{\mathcal{H}}}{\partial x}(t) \left(
X(t)-\hat{X}(t)\right) \leq 0,
\end{equation}
we use the Arrow condition \eqref{eq3.8A} and the maximum condition \eqref{eq3.9}, proceeding as in the proof of Corollary 11.3 in \cite{OS1}. We omit the details.
%&  \leq\mathbb{E}\left[  \int_{0}^{T} \frac{\partial \hat{\mathcal{H}}}{\partial v}(t) \left(  u(t)-\hat{u}(t)\right)  dt\right] \\
%&  =\mathbb{E}\left[  \int_{0}^{T}\mathbb{E}\left[  \left.   \frac{\partial \hat{\mathcal{H}}}{\partial v}(t) \left(  u(t)-\hat{u}(t)\right) \right\vert \mathcal{G}_{t}\right]  dt\right] \nonumber\\
%&  =\mathbb{E}\left[  \int_{0}^{T}\mathbb{E}\left[  \left.   \frac{\partial \hat{\mathcal{H}}}{\partial v}(t) \right\vert \mathcal{G}
%_{t}\right]  \left(  u(t)-\hat{u}(t)\right)  dt\right]  \leq0\text{ by
%}\left(  \ref{1.7}\right)  .
%\end{align*}
\fproof

\section{A Necessary Maximum Principle}

The sufficient maximum principle proved in the previous section has the drawback that the required concavity conditions are not always satisfied. It is therefore important also to have a maximum principle which does not need this assumption. In the following result, a \emph{necessary} maximum principle, the concavity conditions are replaced by conditions related to the space of admissible controls and the existence of the derivative process. The details are as follows:

For each given $t \in [0,T]$ let $\alpha=\alpha_t$ be a bounded $\mathcal{G}_t$ - measurable random variable, let $h \in [T-t,T]$ and define

\begin{equation}\label{eq4.1}
\beta(s):=\alpha1_{\left[  t,t+h\right]  }(s);s\in\left[  0,T\right].
%\label{1.12}%
\end{equation}

Assume that
\begin{equation}
u + \lambda \beta \in \mathcal{A}
\end{equation}
for all such  $\alpha$ and all $u\in \mathcal{A}$, and all non-zero $\lambda$ sufficiently small. Assume that the \emph{derivative process} $Y(t)$, defined by
\begin{equation}
Y(t)=\left.  \frac{d}{d\lambda}X^{(u+\lambda\beta)}(t)\right\vert _{\lambda
=0}, \label{1.13}%
\end{equation}
exists.

Then we see that%
\begin{align}
Y(t) &  =%
%TCIMACRO{\dint \limits_{0}^{t}}%
%BeginExpansion
{\displaystyle\int\limits_{0}^{t}}
%EndExpansion
\left(  \frac{\partial b}{\partial x}(t,s)Y(s)+\frac{\partial b}{\partial
u}(t,s)\beta(s)\right)  ds\nonumber\\
&  +%
%TCIMACRO{\dint \limits_{0}^{t}}%
%BeginExpansion
{\displaystyle\int\limits_{0}^{t}}
%EndExpansion
\left(  \frac{\partial\sigma}{\partial x}(t,s)Y(s)+\frac{\partial\sigma
}{\partial u}(t,s)\beta(s)\right)  dB(s)\nonumber\\
&  +%
%TCIMACRO{\dint \limits_{0}^{t}}%
%BeginExpansion
{\displaystyle\int\limits_{0}^{t}}
%EndExpansion%
%TCIMACRO{\dint \limits_{\mathbb{R}}}%
%BeginExpansion
{\displaystyle\int\limits_{\mathbb{R}}}
%EndExpansion
\left(  \frac{\partial\gamma}{\partial x}(t,s,\zeta)Y(s)+\frac{\partial\gamma
}{\partial u}(t,s,\zeta)\beta(s)\right)  \tilde{N}(ds,d\zeta),
\end{align}

and hence%

\begin{align}
dY(t) &  =\left[  \frac{\partial b}{\partial x}(t,t)Y(t)+\frac{\partial
b}{\partial u}(t,t)\beta(t)+%
%TCIMACRO{\dint \limits_{0}^{t}}%
%BeginExpansion
{\displaystyle\int\limits_{0}^{t}}
%EndExpansion
\left(  \frac{\partial^{2}b}{\partial t\partial x}(t,s)Y(s)+\frac{\partial
^{2}b}{\partial t\partial u}(t,s)\beta(s)\right)  ds\right.  \nonumber\\
&  +%
%TCIMACRO{\dint \limits_{0}^{t}}%
%BeginExpansion
{\displaystyle\int\limits_{0}^{t}}
%EndExpansion
\left(  \frac{\partial^{2}\sigma}{\partial t\partial x}(t,s)Y(s)+\frac
{\partial^{2}\sigma}{\partial t\partial u}(t,s)\beta(s)\right)
dB(s)\nonumber\\
&  \left.  +%
%TCIMACRO{\dint \limits_{0}^{t}}%
%BeginExpansion
{\displaystyle\int\limits_{0}^{t}}
%EndExpansion%
%TCIMACRO{\dint \limits_{\mathbb{R}}}%
%BeginExpansion
{\displaystyle\int\limits_{\mathbb{R}}}
%EndExpansion
\left(  \frac{\partial^{2}\gamma}{\partial t\partial x}(t,s,\zeta
)Y(s)+\frac{\partial^{2}\gamma}{\partial t\partial u}(t,s,\zeta)\beta
(s)\right)  \tilde{N}(ds,d\zeta)\right]  dt\nonumber\\
&  +\left(  \frac{\partial\sigma}{\partial x}(t,t)Y(t)+\frac{\partial\sigma
}{\partial u}(t,t)\beta(t)\right)  dB(t)\nonumber\\
&  +%
%TCIMACRO{\dint \limits_{\mathbb{R}}}%
%BeginExpansion
{\displaystyle\int\limits_{\mathbb{R}}}
%EndExpansion
\left(  \frac{\partial\gamma}{\partial x}(t,t,\zeta)Y(t)+\frac{\partial\gamma
}{\partial u}(t,t,\zeta)\beta(t)\right)  \tilde{N}(dt,d\zeta).\label{1.14}%
\end{align}
We are now ready to formulate the result

\begin{theorem} [Necessary maximum principle]

Suppose that $\hat{u}\in$ $\mathcal{A}$ is such that, for all $\beta$ as in \eqref{eq4.1},
\begin{equation}\label{eq4.5}
\left.  \frac{d}{d\lambda}J(\hat{u}+\lambda\beta)\right\vert _{\lambda=0}=0
\end{equation}
and the corresponding solution $\hat{X}(t),(\hat{p}(t),\hat{q}(t),\hat{r}(t,\cdot))$ of \eqref{eq1.1}
 and \eqref{eq3.4} exists.

Then,
\begin{equation}\label{eq4.6}
\mathbb{E}\left[  \left.  \dfrac{\partial\mathcal{H}}{\partial u}(t)\right\vert \mathcal{G}_{t}\right] _{u=u(t)} =0.
\end{equation}
Conversely, if \eqref{eq4.6} holds, then \eqref{eq4.5} holds.
\end{theorem}

\dproof

For simplicity of notation we drop the "hat" everywhere and write $u$ in stead of $\hat{u}$, $X$ in stead of $\hat{X}$ etc in the following.\\
By considering a suitable increasing family of stopping times converging to
$T$, we may assume that all the local martingales appearing in the proof below
are martingales. See the proof of Theorem 2.1 in \cite{OS2} for details.\\
Now consider
\begin{equation}\label{1.15}%
\begin{array}
[c]{l}%
\left.  \dfrac{d}{d\lambda}J(u+\lambda\beta)\right\vert _{\lambda=0}\\
=\mathbb{E}\left[
%TCIMACRO{\dint \limits_{0}^{T}}%
%BeginExpansion
{\displaystyle\int\limits_{0}^{T}}
%EndExpansion
\left\{  \dfrac{\partial f}{\partial x}(t)Y(t)+\dfrac{\partial f}{\partial
u}(t)\beta(t)\right\}  dt+g^{\prime}(X(T))Y(T)\right].
\end{array}
\end{equation}%

Applying the It\^{o} Formula,
%to $g^{\prime }(X(T)Y(T))$,
we get
\begin{align*}
&
\begin{array}
[c]{l}%
\mathbb{E}\left[  g^{\prime}(X(T))Y(T)\right]  =\mathbb{E}\left[
p(T)Y(T)\right]  \\
=\mathbb{E}\left[
%TCIMACRO{\dint \limits_{0}^{T}}%
%BeginExpansion
{\displaystyle\int\limits_{0}^{T}}
%EndExpansion
p(t)\left(  \dfrac{\partial b}{\partial x}(t,t)Y(t)+\dfrac{\partial
b}{\partial u}(t,t)\beta(t)\right)  dt\right.  \\
+%
%TCIMACRO{\dint \limits_{0}^{T}}%
%BeginExpansion
{\displaystyle\int\limits_{0}^{T}}
%EndExpansion
p(t)\left\{
%TCIMACRO{\dint \limits_{0}^{t}}%
%BeginExpansion
{\displaystyle\int\limits_{0}^{t}}
%EndExpansion
\left(  \frac{\partial^{2}b}{\partial t\partial x}(t,s)Y(s)+\frac{\partial
^{2}b}{\partial t\partial u}(t,s)\beta(s)\right)  ds\right\}  dt\\
+%
%TCIMACRO{\dint \limits_{0}^{T}}%
%BeginExpansion
{\displaystyle\int\limits_{0}^{T}}
%EndExpansion
p(t)\left\{
%TCIMACRO{\dint \limits_{0}^{t}}%
%BeginExpansion
{\displaystyle\int\limits_{0}^{t}}
%EndExpansion
\left(  \frac{\partial^{2}\sigma}{\partial t\partial x}(t,s)Y(s)+\frac
{\partial^{2}\sigma}{\partial t\partial u}(t,s)\beta(s)\right)  dB(s)\right\}
dt\\
+%
%TCIMACRO{\dint \limits_{0}^{T}}%
%BeginExpansion
{\displaystyle\int\limits_{0}^{T}}
%EndExpansion
p(t)\left\{
%TCIMACRO{\dint \limits_{0}^{t}}%
%BeginExpansion
{\displaystyle\int\limits_{0}^{t}}
%EndExpansion%
%TCIMACRO{\dint \limits_{\mathbb{R}}}%
%BeginExpansion
{\displaystyle\int\limits_{\mathbb{R}}}
%EndExpansion
\left(  \frac{\partial^{2}\gamma}{\partial t\partial x}(t,s,\zeta
)Y(s)+\frac{\partial^{2}\gamma}{\partial t\partial u}(t,s,\zeta)\beta
(s)\right)  \tilde{N}(ds,d\zeta)\right\}  dt\\
-%
%TCIMACRO{\dint \limits_{0}^{T}}%
%BeginExpansion
{\displaystyle\int\limits_{0}^{T}}
%EndExpansion
Y(t)
%\mathbb{E}\left[  \left.
\frac{\partial\mathcal{H}}{\partial x}(t)
%\right\vert \mathcal{F}_{t}\right]
dt+%
%TCIMACRO{\dint \limits_{0}^{T}}%
%BeginExpansion
{\displaystyle\int\limits_{0}^{T}}
%EndExpansion
q(t)\left(  \frac{\partial\sigma}{\partial x}(t,t)Y(t)+\frac{\partial\sigma
}{\partial u}(t,t)\beta(t)\right)  dt
\end{array}
\\
&  \left.  +%
%TCIMACRO{\dint \limits_{0}^{T}}%
%BeginExpansion
{\displaystyle\int\limits_{0}^{T}}
%EndExpansion%
%TCIMACRO{\dint \limits_{\mathbb{R}}}%
%BeginExpansion
{\displaystyle\int\limits_{\mathbb{R}}}
%EndExpansion
r(t,\zeta)\left(  \frac{\partial\gamma}{\partial x}(t,t,\zeta)Y(t)+\frac
{\partial\gamma}{\partial u}(t,t,\zeta)\beta(t)\right)  \nu(d\zeta)dt\right].
\end{align*}

From \eqref{eq2.11}, \eqref{eq2.12} and \eqref{eq2.13}, we have%
\begin{align*}
&
\begin{array}
[c]{l}%
\mathbb{E}\left[  p(T)Y(T)\right]  \\
=\mathbb{E}\left[
%TCIMACRO{\dint \limits_{0}^{T}}%
%BeginExpansion
{\displaystyle\int\limits_{0}^{T}}
%EndExpansion
\left\{  \dfrac{\partial b}{\partial x}(t,t)p(t)+%
%TCIMACRO{\dint \limits_{t}^{T}}%
%BeginExpansion
{\displaystyle\int\limits_{t}^{T}}
%EndExpansion
\left(  \frac{\partial^{2}b}{\partial s\partial x}(s,t)p(s)+\frac{\partial
^{2}\sigma}{\partial s\partial x}(s,t)\mathbb{E}[D_{t}p(s)\mid\mathcal{F}_t]\right.  \right.  \right.  \\
\left.  \left.  +%
%TCIMACRO{\dint \limits_{\mathbb{R}}}%
%BeginExpansion
{\displaystyle\int\limits_{\mathbb{R}}}
%EndExpansion
\frac{\partial^{2}\gamma}{\partial s\partial x}(s,t,\zeta)\mathbb{E}[D_{t,\zeta}%
p(s)\mid\mathcal{F}_t]\nu(d\zeta)\right)  ds\right\}  Y(t)dt\\
\mathbb{+}%
%TCIMACRO{\dint \limits_{0}^{T}}%
%BeginExpansion
{\displaystyle\int\limits_{0}^{T}}
%EndExpansion
\left\{  \dfrac{\partial b}{\partial u}(t,t)p(t)+%
%TCIMACRO{\dint \limits_{t}^{T}}%
%BeginExpansion
{\displaystyle\int\limits_{t}^{T}}
%EndExpansion
\left(  \frac{\partial^{2}b}{\partial s\partial u}(s,t)p(s)+\frac{\partial
^{2}\sigma}{\partial s\partial u}(s,t)\mathbb{E}[D_{t}p(s)\mid\mathcal{F}_t]\right.  \right.  \\
\left.  \left.  +%
%TCIMACRO{\dint \limits_{\mathbb{R}}}%
%BeginExpansion
{\displaystyle\int\limits_{\mathbb{R}}}
%EndExpansion
\frac{\partial^{2}\gamma}{\partial s\partial u}(s,t,\zeta)\mathbb{E}[D_{t,\zeta}%
p(s)\mid\mathcal{F}_t]\nu(d\zeta)\right)  ds\right\}  \beta(t)dt\\
-%
%TCIMACRO{\dint \limits_{0}^{T}}%
%BeginExpansion
{\displaystyle\int\limits_{0}^{T}}
%EndExpansion
\frac{\partial\mathcal{H}}{\partial x}(t)Y(t)dt+%
%TCIMACRO{\dint \limits_{0}^{T}}%
%BeginExpansion
{\displaystyle\int\limits_{0}^{T}}
%EndExpansion
\left(  \frac{\partial\sigma}{\partial x}(t,t)Y(t)+\frac{\partial\sigma
}{\partial u}(t,t)\beta(t)\right)  q(t)dt
\end{array}
\\
&  \left.  +%
%TCIMACRO{\dint \limits_{0}^{T}}%
%BeginExpansion
{\displaystyle\int\limits_{0}^{T}}
%EndExpansion%
%TCIMACRO{\dint \limits_{\mathbb{R}}}%
%BeginExpansion
{\displaystyle\int\limits_{\mathbb{R}}}
%EndExpansion
\left(  \frac{\partial\gamma}{\partial x}(t,t,\zeta)Y(t)+\frac{\partial\gamma
}{\partial u}(t,t,\zeta)\beta(t)\right)  r(t,\zeta)\nu(d\zeta)dt\right].
\end{align*}

Using the definition of $\mathcal{H}$ in \eqref{eq3.3} and the definition of $\beta$, we obtain%

\begin{equation} \label{eq4.8}
\left.  \dfrac{d}{d\lambda}J(u+\lambda\beta)\right\vert _{\lambda
=0}=\mathbb{E}\left[
{\displaystyle\int\limits_{0}^{T}}
\frac{\partial\mathcal{H}}{\partial u}(s)\beta(s)ds\right]
 =\mathbb{E}\left[
{\displaystyle\int\limits_{t}^{t+h}}
\frac{\partial\mathcal{H}}{\partial u}(s)ds\alpha\right].
\end{equation}

Now suppose that
\begin{equation} \label{eq4.9}
\left.  \dfrac{d}{d\lambda}J(u+\lambda\beta)\right\vert _{\lambda
=0}=0.
\end{equation}

Differentiating the right-hand side of \eqref{eq4.8} at $h=0$, we get%
\[
\mathbb{E}\left[  \frac{\partial\mathcal{H}}{\partial u}(t)\alpha\right]  =0.
\]

Since this holds for all bounded $\mathcal{G}_{t}$-measurable $\alpha$, we have%
\begin{equation} \label{eq4.10}
\mathbb{E}\left[  \left.  \frac{\partial\mathcal{H}}{\partial u}(t)\right\vert
\mathcal{G}_{t}\right]  =0.
\end{equation}

Conversely, if we assume that \eqref{eq4.10} holds, then we obtain \eqref{eq4.9} by reversing the argument we used to obtain \eqref{eq4.8}.

\fproof

\section{Applications}

\subsection{The Case when the Coefficients do not Depend on $x$}

Consider the case when the coefficients do not depend on $x,$ i.e., the system
has the form:
\begin{eqnarray}\label{eq5.1}
X(t) &=&\xi (t)+\int_{0}^{t}b(t,s,u(s))ds+\int_{0}^{t}\sigma (t,s,u(s))dB(s)
\label{4.1} \\
&&+\int_{0}^{t}\int_{
\mathbb{R}
}\gamma (t,s,u(s),\zeta )\tilde{N}(ds,d\zeta )  \nonumber
\end{eqnarray}

\noindent
with performance functional

\begin{equation}
J(u)=\mathbb{E}\left[ \int_{0}^{T}f(t,u(t))dt+g(X(T))\right].   \label{4.2}
\end{equation}
This special case is of interest because any linear SDE with delay can be written on this form. See \cite{OZ3} and the references therein.
$\bigskip $\\
In this case the Hamiltonian $\mathcal{H}$ given in \eqref{eq3.3} takes
the form
\begin{eqnarray} \label{eq4.3}
&&\mathcal{H}(t,v,p(\cdot ),q(\cdot ),r(\cdot )) \nonumber\\
&=&f(t,v)+b(t,t,v)p(t)+\sigma (t,t,v)q(t)+\int_{
\mathbb{R}
}\gamma (t,t,v,\zeta )r(t,\zeta )\nu (d\zeta ) \nonumber\\
&&+\int_{t}^{T}\frac{\partial b}{\partial s}(s,t,v)p(s)ds+\int_{t}^{T}\frac{
\partial \sigma }{\partial s}(s,t,v)\mathbb{E}[D_{t}p(s)\mid\mathcal{F}_t]ds \nonumber\\
&&+\int_{t}^{T}\int_{
\mathbb{R}
}\frac{\partial \gamma }{\partial s}(s,t,v,\zeta )\mathbb{E}[D_{t,\zeta }
p(s)\mid\mathcal{F}_t]\nu (d\zeta )ds.
\end{eqnarray}

The BSDE \eqref{eq3.4} for the adjoint variables $p,q,r$ gets the form

\begin{equation}
\left\{
\begin{array}{l}
dp(t)=q(t)dB(t)+\int_{
\mathbb{R}
}r(t,\zeta )\tilde{N}(dt,d\zeta );0\leq t\leq T \\
p(t)=g^{\prime }(X(T)),
\end{array}
\right.   \label{4.4}
\end{equation}

which has the solution
\begin{eqnarray}
p(t) &=&\mathbb{E}[g^{\prime }(X(T))\mid \mathcal{F}_{t}]  \label{4.5} \\
q(t) &=&D_{t}p(t)=\mathbb{E}[D_{t}g^{\prime }(X(T))\mid
\mathcal{F}_{t}]  \label{4.6} \\
r(t,\zeta ) &=&D_{t,\zeta }p(t)=\mathbb{E}[D_{t,\zeta }
g^{\prime }(X(T))\mid \mathcal{F}_{t}].  \label{4.7}
\end{eqnarray}

Substituting (\ref{4.5})-(\ref{4.7}) into \eqref{eq4.3} we get
\[
\mathbb{E}[\mathcal{H}(t,v,p(\cdot ),q(\cdot ),r(\cdot ))\mid \mathcal{F}
_{t}]=\mathbb{E}[\mathcal{H}_{0}(t,v,p,q,r)\mid \mathcal{F}_{t}],
\]

where
\begin{eqnarray}
\mathcal{H}_{0}(t,v,p,q,r) &=&f(t,v)+b(t,t,v)g^{\prime }(X(T))+\sigma
(t,t,v)\mathbb{E}[D_{t}g^{\prime }(X(T))\mid \mathcal{F}_t] \label{4.9} \nonumber\\
&&+\int_{
\mathbb{R}
}\gamma (t,t,v,\zeta )\mathbb{E}[D_{t,\zeta }g^{\prime }(X(T))\mid \mathcal{F}_t] \nu (d\zeta )
\nonumber \\
&&+\int_{t}^{T}\frac{\partial b}{\partial s}(s,t,v)g^{\prime
}(X(T))ds+\int_{t}^{T}\frac{\partial \sigma }{\partial s}(s,t,v)\mathbb{E}[D_{t}g^{\prime }(X(T))\mid \mathcal{F}_t]ds  \nonumber \\
&&+\int_{t}^{T}\int_{
\mathbb{R}
}\frac{\partial \gamma }{\partial s}(s,t,v,\zeta )\mathbb{E}[D_{t,\zeta }g^{\prime }(X(T))\mid \mathcal{F}_t]\nu (d\zeta )ds.
\end{eqnarray}

Performing the ds-integrals we see that $\mathcal{H}_{0}(t,v,p,q,r)$\
reduces to
\begin{eqnarray}
\mathbb{H}_{0}(t,v,X(T)) &:&=f(t,v)+b(T,t,v)g^{\prime }(X(T))  \nonumber \\
&&+\sigma (T,t,v)\mathbb{E}[D_{t}g^{\prime }(X(T))\mid \mathcal{F}_t] \nonumber\\
&&+\int_{
\mathbb{R}
}\gamma (T,t,v,\zeta )\mathbb{E}[D_{t,\zeta }g^{\prime }(X(T))\mid \mathcal{F}_t]\nu (d\zeta ).
\end{eqnarray}

We conclude that, in this case, we have the following maximum principles:

\begin{theorem}[Sufficient maximum principle II]

Suppose that the coefficients $f(t,v),b(t,s,v),\sigma (t,s,v)$ and $\gamma
(t,s,v,\zeta )$ of the stochastic control system (\ref{4.1})-(\ref{4.2}) do
not depend on $x.$

Let $\hat{u}\in \mathcal{A}$ with associated solution $\hat{X}$%
\ of (\ref{4.1}). Suppose that the functions

\begin{align}
x\rightarrow g(x)
\end{align}
$and$
\begin{align}
v\rightarrow \mathbb{H}_{0}(t,v,\hat{X}(T))
\end{align}
are concave and that, for all $t,$
\begin{equation}
\max_{v\in \mathbb{U}}\mathbb{E}\left[ \mathbb{H}_{0}(t,v,\hat{X}(T))\mid
\mathcal{G}_{t}\right] =\mathbb{E}\left[ \mathbb{H}_{0}(t,\hat{u}(t),\hat{X}
(T))\mid \mathcal{G}_{t}\right] .  \label{4.12}
\end{equation}

Then, $\hat{u}$\ is an optimal control, i.e.,
\begin{equation}
\sup_{u\in \mathcal{A}}J(u)=J(\hat{u}).  \label{4.13}
\end{equation}
\end{theorem}

\begin{theorem}[Necessary maximum principle II]

Let $X(t)$ and $J(u)$ be as in Theorem $5.1.$

Let $\hat{u}\in \mathcal{A}$ with associated solution $\hat{X}$
of \eqref{eq5.1}.

Then the following, (i) and (ii), are equivalent:

(i) $\hat{u}$\ is a critical point for $J(u)$,i.e.,
\[
\left. \frac{d}{dy}J(\hat{u}+yw)\right\vert _{y=0}=0
\]

for all processes $w$ such that $\hat{u}+yw\in \mathcal{A}$ for all $y$ small enough.

(ii)
\[
\mathbb{E}\left[ \frac{\partial \mathbb{H}_{0}}{\partial v}(t,v,\hat{X}
(T))\mid \mathcal{G}_{t}\right] _{v=\hat{u}(t)}=0.
\]
\end{theorem}

\begin{remark}
Theorem 5.2 is identical to Theorem 3.2 in \cite{OZ3}. However, the method in \cite{OZ3} is different, being based on perturbation techniques and
complicated stochastic expansions. In the general case the necessary maximum
principle of \cite{OZ3} is completely different from our Theorem 3.1. There is no
corresponding \emph{sufficient} maximum principle in \cite{OZ3}.
\end{remark}

\subsection{Optimal Investment in a Financial Market Modeled by a Volterra
Equation}

Consider a financial market with the following two investment possibilities:

(i) A risk free asset with unit price $S_{0}(t)=1$; $t\geq 0$

(ii) A risky asset, in which investments have long term (memory) effects, in
the following sense:

If we at time $s\geq 0$ decide to invest the fraction $\pi (s)$ of the
current total wealth $X(s)$ in this asset, then we assume that the wealth $%
X(t)=X_{\pi }(t)$ at time $t$ is described by the linear stochastic Volterra
equation
\begin{equation}
X(t)=x+\int_{0}^{t}b_{0}(t,s)\pi (s)X(s)ds+\int_{0}^{t}\sigma _{0}(t,s)\pi
(s)X(s)dB(s);t\geq 0  \label{eq5.12}
\end{equation}

or, in differential form,

\begin{equation}
\left\{
\begin{array}{l}
dX(t)=b_{0}(t,t)\pi (t)X(t)dt+\sigma _{0}(t,t)\pi (t)X(t)dB(t) \\
+\left[ \int_{0}^{t}\frac{\partial b_{0}}{\partial t}(t,s)\pi
(s)X(s)ds+\int_{0}^{t}\frac{\partial \sigma _{0}}{\partial t}(t,s)\pi
(s)X(s)dB(s)\right] dt;t\geq 0 \\
X(0)=x.
\end{array}%
\right.  \label{eq5.13}
\end{equation}

Thus we see that \eqref{eq5.13} differs from the classical
Black-Scholes type of wealth equation by the last two integral terms on the
right hand side. These terms represent long term (memory) effects of the
investment strategy $\pi (\cdot )$.

We assume that $b_{0}(t,s)=b_{0}(t,s,\omega )$ and $\sigma _{0}(t,s)=\sigma
_{0}(t,s,\omega )$ are given bounded processes, and that $b_{0}(t,s)$ and $%
\sigma _{0}(t,s)$ are $\mathcal{F}_{s}$-measurable for all $s,t$ and
differentiable with bounded derivatives with respect to $t$ for all $s$, a.s.

We also assume that \begin{equation}\label{eq5.15a}
\sigma_0(t,s) \geq c_0  \text{ a.s. for all } t,s \in[0,T]\\
\text{ for some constant } c_0 >0.
\end{equation}
We choose $\mathbb{G}=\mathbb{F}$ in this example and we say that $\pi $ is
admissible and write $\pi \in \mathcal{A}$ if $\pi $ is $\mathbb{F}$%
-adapted, $\pi \in L^{2}(d\lambda \times dP)$ and equation \eqref{eq5.12} has a unique solution with $\pi X\in L^{2}(d\lambda \times dP).$

We assume that $x>0.$ If $\pi \in \mathcal{A}$, then it follows that $X_{\pi }(t)>0
$ for all $t\in \lbrack 0,T].$ To see this, note that from \eqref{eq5.13} we get%
\begin{equation}
\begin{array}{l}
X_{\pi }(t)=x\exp \left( \int_{0}^{t}\sigma _{0}(s,s)\pi (s)dB(s)\right.  \\
\left. +\int_{0}^{t}\left\{ b_{0}(s,s)\pi (s)-\frac{1}{2}\sigma
_{0}^{2}(s,s)\pi ^{2}(s)+\alpha (s)\right\} ds\right) >0,%
\end{array}%
\end{equation}

where
\begin{equation*}
\alpha (s):=\int_{0}^{s}\frac{\partial b_{0}}{\partial s}(s,r)\pi
(r)X(r)dr+\int_{0}^{s}\frac{\partial \sigma _{0}}{\partial s}(s,r)\pi
(r)X(r)dB(r).
\end{equation*}
\noindent
We now study the following optimal investment problem:

Find $\hat{\pi}\in \mathcal{A}$ such that
\begin{equation}\label{eq5.18}
\underset{\pi \in \mathcal{A}}{\sup }\mathbb{E}\left[ U(X_{\pi }(T))\right] =%
\mathbb{E}\left[ U(X_{\hat{\pi}}(T))\right],
\end{equation}
\noindent
where $U:\left[ 0,\infty \right) \rightarrow $ $\left[-\infty ,\infty
\right)$ is a given \textit{utility }function, assumed to be strictly
increasing, concave, and $C^{1}$ on $(0,\infty)$.
\noindent
This is a control problem of the type studied in Section $3$ and $4$, and we
apply the results from there:

The Hamiltonian $\mathcal{H}$ given by \eqref{eq3.3} gets the form%
\begin{equation}
\begin{array}{l}
\mathcal{H}(t,x,\pi,p,q)=b_{0}(t,t)\pi xp+\sigma _{0}(t,t)\pi xq \\
+\int_{t}^{T}\frac{\partial b_{0}}{\partial s}(s,t)\pi xp(s)ds+\int_{t}^{T}%
\frac{\partial \sigma _{0}}{\partial s}(s,t)\pi x\mathbb{E}[D_{t}p(s)\mid \mathcal{F}_t]ds%
\end{array}
\label{eq5.16}
\end{equation}

Suppose there exists an optimal control $\hat{\pi}\in \mathcal{A}$ for \eqref{eq5.18} with corresponding $\hat{X},\hat{p},\hat{q}.$ Then,
\begin{equation*}
\mathbb{E}\left[ \left. \frac{\partial }{\partial \pi }\mathcal{H}(t,\hat{X}%
(t),\pi ,\hat{p},\hat{q})\right\vert \mathcal{F}_{t}\right] _{\pi =\hat{\pi}%
(t)}=0
\end{equation*}

i.e.,%
\begin{equation*}
\begin{array}{l}
\mathbb{E}\left[ b_{0}(t,t)\hat{X}(t)\hat{p}(t)+\sigma _{0}(t,t)\hat{X}(t)%
\hat{q}(t)\right.  \\
\left. \left. +\int_{t}^{T}\frac{\partial b_{0}}{\partial s}(s,t)\hat{X}(t)%
\hat{p}(s)ds+\int_{t}^{T}\frac{\partial \sigma _{0}}{\partial s}(s,t)\hat{X}%
(t)\mathbb{E}[D_{t}\hat{p}(s)\mid\mathcal{F}_t]ds\right\vert \mathcal{F}_{t}\right] =0.%
\end{array}%
\end{equation*}

Since $\hat{X}(t)>0$, this is equivalent to%
\begin{equation}\label{eq5.17}
\begin{array}{l}
b_{0}(t,t)\hat{p}(t)+\sigma _{0}(t,t)\hat{q}(t) \\
+\mathbb{E}\left[ \left. \int_{t}^{T}\left\{ \frac{\partial b_{0}}{\partial s%
}(s,t)\hat{p}(s)+\int_{t}^{T}\frac{\partial \sigma _{0}}{\partial s}%
(s,t)\mathbb{E}[D_{t}\hat{p}(s)\mid\mathcal{F}_t]\right\} ds\right\vert \mathcal{F}_{t}\right] =0.%
\end{array}%
\end{equation}

We deduce that the corresponding BSDE (3.4) reduces to%
\begin{equation}
\left\{
\begin{array}{l}
d\hat{p}(t)=\hat{q}(t)dB(t);0\leq t\leq T \\
\hat{p}(T)=U^{\prime }(\hat{X}(T)),%
\end{array}%
\right.
\end{equation}

which has the unique solution%
\begin{equation}
\hat{p}(t)=\mathbb{E}\left[ \left. U^{\prime }(\hat{X}(T))\right\vert
\mathcal{F}_{t}\right] ,\hat{q}(t)=D_{t}\hat{p}(t).
\end{equation}

Substituted into \eqref{eq5.17}, this gives the equation%
\begin{equation} \label{eq5.20}
\begin{array}{l}
\mathbb{E}\left[ b_{0}(t,t)U^{\prime }(\hat{X}(T))+\sigma
_{0}(t,t)D_{t}U^{\prime }(\hat{X}(T))\right.  \\
+\int_{t}^{T}\frac{\partial b_{0}}{\partial s}(s,t)\mathbb{E}\left[ \left.
U^{\prime }(\hat{X}(T))\right\vert \mathcal{F}_{s}\right] ds \\
\left. \left. +\int_{t}^{T}\frac{\partial \sigma _{0}}{\partial s}(s,t)%
\mathbb{E}\left[ \left. D_{t}U^{\prime }(\hat{X}(T))\right\vert \mathcal{F}%
_{s}\right] ds\right\vert \mathcal{F}_{t}\right] =0,
\end{array}%
\end{equation}
\noindent
where we have used that%
\begin{equation}
D_{t}\mathbb{E}\left[ \left. U^{\prime }(\hat{X}(T))\right\vert \mathcal{F}%
_{t}\right] =\mathbb{E}\left[ \left. D_{t}U^{\prime }(\hat{X}(T))\right\vert
\mathcal{F}_{t}\right] ,
\end{equation}
which is an identity that follows easily from the definition \eqref{eq2.5a} of the Malliavin derivative.
Equation \eqref{eq5.20} can be simplified to%
\begin{equation}
\begin{array}{l}
b_{0}(t,t)\mathbb{E}\left[ \left. U^{\prime }(\hat{X}(T))\right\vert
\mathcal{F}_{t}\right] +\sigma _{0}(t,t)\mathbb{E}\left[ \left.
D_{t}U^{\prime }(\hat{X}(T))\right\vert \mathcal{F}_{t}\right]  \\
+\mathbb{E}\left[ \left. \int_{t}^{T}\frac{\partial b_{0}}{\partial s}%
(s,t)U^{\prime }(\hat{X}(T))ds\right\vert \mathcal{F}_{t}\right]  \\
+\mathbb{E}\left[ \left. \int_{t}^{T}\frac{\partial \sigma _{0}}{\partial s}%
(s,t)D_{t}U^{\prime }(\hat{X}(T))ds\right\vert \mathcal{F}_{t}\right] =0,%
\end{array}%
\end{equation}

or%
\begin{equation}\label{eq5.23}
\sigma _{0}(T,t)D_{t}\mathbb{E}\left[ \left. U^{\prime }(\hat{X}%
(T))\right\vert \mathcal{F}_{t}\right] +b_{0}(T,t)\mathbb{E}\left[ \left.
U^{\prime }(\hat{X}(T))\right\vert \mathcal{F}_{t}\right] =0.
\end{equation}

By \eqref{eq5.15a} we see that \eqref{eq5.23} can be written

\begin{equation}
\frac{D_{t}Y(t)}{Y(t)}=-\frac{b_{0}(T,t)}{\sigma _{0}(T,t)},  \label{eq5.25}
\end{equation}

where%
\begin{equation}
Y(t)=\mathbb{E}\left[ \left. U^{\prime }(\hat{X}(T))\right\vert \mathcal{F}%
_{t}\right].  \label{eq5.26}
\end{equation}

By the chain rule for Malliavin derivatives, we deduce from \eqref{eq5.26}
that%
\begin{equation}
D_{t}(\ln Y(t))=-\frac{b_{0}(T,t)}{\sigma _{0}(T,t)}.  \label{eq5.27}
\end{equation}

On the other hand, since $Y(t)$\ is a positive martingale, there exists an
adapted process $\theta _{0}(t)$ such that

\begin{equation*}
dY(t)=\theta _{0}(t)Y(t)dB(t)
\end{equation*}

i.e.,
\begin{equation}
Y(t)=Y(0)\exp (\int_{0}^{t}\theta _{0}(s)dB(s)-\frac{1}{2}\int_{0}^{t}\theta
_{0}^{2}(s)ds). \label{eq5.28}
\end{equation}

From \eqref{eq5.28} we get
\begin{equation}
D_{t}(\ln Y(t))=D_{t}(\int_{0}^{t}\theta _{0}(s)dB(s)-\frac{1}{2}%
\int_{0}^{t}\theta _{0}^{2}(s)ds)=\theta _{0}(t),  \label{eq5.29}
\end{equation}

since%
\begin{equation*}
D_{t}\theta _{0}(s)=D_{t}(\theta _{0}^{2}(s))=0\text{ for all }s<t\text{
(because }\theta _{0}\ \text{is adapted).}
\end{equation*}

Comparing \eqref{eq5.27} and \eqref{eq5.29} we conclude that
\begin{equation}
\theta _{0}(t)=-\frac{b_{0}(T,t)}{\sigma _{0}(T,t)}  \label{eq5.30}
\end{equation}

and hence, by \eqref{eq5.26},

\begin{equation}
\mathbb{E}\left[ \left. U^{\prime }(\hat{X}(T))\right\vert \mathcal{F}_{t}%
\right] =Y(t)=\mathbb{E}\left[ U^{\prime }(\hat{X}(T))\right] \exp
(\int_{0}^{t}\theta _{0}(s)dB(s)-\int_0^t \frac{1}{2}\theta _{0}^{2}(s)ds).
\label{eq5.31}
\end{equation}

It remains to find the constant
\begin{equation}
c:=\mathbb{E}\left[ U^{\prime }(\hat{X}(T))\right] .  \label{eq5.32}
\end{equation}

From \eqref{eq5.31} with $t=T$ we get

\begin{equation}
\hat{X}(T)=(U^{\prime })^{-1}(c\exp (\int_{0}^{T}\theta _{0}(s)dB(s)-\frac{1%
}{2}\int_{0}^{T}\theta _{0}^{2}(s)ds))=\colon F(c).  \label{eq5.33}
\end{equation}

On the other hand, if we define

\begin{equation}
\hat{Z}_{c}(t,s):=\sigma _{0}(t,s)\hat{\pi}(s)\hat{X}(s),  \label{eq5.34}
\end{equation}
\noindent
then by \eqref{eq5.12}, the pair $(\hat{X},\hat{Z}_{c})$ solves
the following \emph{(Yong type) backward stochastic Volterra integral equation (BSVIE)}

\begin{equation}\label{eq5.35}
\hat{X}(t)=F(c)-\int_{t}^{T}\frac{b_{0}(t,s)\hat{Z}_{c}(t,s)}{\sigma
_{0}(t,s)}ds-\int_{t}^{T}\hat{Z}_{c}(t,s)dB(s);0\leq s\leq T.
\end{equation}

By Theorem 3.2 in \cite{Yong1} the solution of this equation is
unique. Putting $t=0$ and taking expectation in \eqref{eq5.35}, we get

\begin{equation}\label{eq5.36}
x=\mathbb{E}[F(c)]-\int_{0}^{T}\mathbb{E}\left[ \frac{b_{0}(t,s)}{\sigma _{0}(t,s)}\hat{Z%
}_{c}(t,s)\right] ds.
\end{equation}

This equation determines implicitly the value of $c$.

Hence by  \eqref{eq5.33} we have found the optimal terminal
wealth $\hat{X}(T).$\ Then, finally we obtain the optimal portfolio $\hat{\pi}
$\ by \eqref{eq5.34}.

Conversely, since the functions $x\rightarrow U(x)$\ and $(x,\pi )\rightarrow
\mathcal{H}(t,x,\pi ,\hat{p},\hat{q})$\ are concave, we see that $\hat{\pi}$%
\ found above satisfies the conditions of Theorem 3.1, and hence $\hat{\pi}$\
is indeed optimal. \\
\newline
We summarize what we have proved as follows:

\begin{theorem}

Assume that $\sigma _{0}(t,s)>0$ is bounded away from $0$, for $s,t\in
\lbrack 0,T].$

Then, the optimal portfolio $\hat{\pi}$\ for the problem \eqref{eq5.18} is%
\begin{equation*}
\hat{\pi}(s)=\frac{\hat{Z}_{c}(t,s)}{\sigma _{0}(t,s)\hat{X}(s)};s\in
\lbrack 0,T],
\end{equation*}
\noindent
where $(\hat{X},\hat{Z}_{c})$ is the unique solution of the BSVIE \eqref{eq5.35}
 with $F$ defined by \eqref{eq5.33}, and the constant $c$ is the solution of \eqref{eq5.36}.
\end{theorem}
\section {Conclusions}

In this paper we study the problem of optimal control of a system described by a stochastic Volterra equation with jumps.
We use Malliavin calculus to obtain both a sufficient and a necessary maximum principle for optimal control of such systems, with partial information. We define a Hamiltonian which involves also the Malliavin derivatives of one of the adjoint processes. This has the advantage that the corresponding adjoint equation becomes (in some way) a standard BSDE, not a Volterra type BSVIE as in \cite{Yong1, TQY} and \cite{SWY}. On the other hand, our BSDE involves the Malliavin derivative of the adjoint process. It is interesting to note that BSDEs involving Malliavin derivatives also appear in connection with optimal control of SDEs with \emph{noisy memory}. See \cite{Dahl}.\\

Our sufficient maximum principle is new, even in the case without jumps. In our general setting also our necessary maximum principle is new. However, in the special case when the coefficients of the state equation do not depend on the state, we show that the necessary maximum principle we obtain, is equivalent to the one in \cite{OZ3}. For more general systems our maximum principle is simpler than the one in \cite{OZ3}. Moreover, there is no sufficient maximum principle in \cite{OZ3}.

In the last part of the paper we illustrate our results by solving an optimal portfolio problem in a financial market with memory, modeled by a stochastic Volterra equation.\\

\paragraph{Acknowledgments}

\begin{itemize}
\item
We are grateful to Nils Christian Framstad for helpful comments.
\item
\thanks{The research leading to these results has
received funding from the European Research Council under the European
Community's Seventh Framework Programme (FP7/2007-2013) / ERC grant agreement
no [228087].}
\item
\thanks{This research was carried out with support of CAS - Centre for Advanced Study, at the Norwegian Academy of Science and Letters, research program SEFE.}
\end{itemize}

\end{document}